\documentclass[11pt,a4paper]{article}

\usepackage{latexsym}
\usepackage{amsmath}
\usepackage{epsfig}
\usepackage{times}
\usepackage{amssymb}
\newtheorem{thr}{\quad Theorem}
\newtheorem{lem}{\quad Lemma}

\newtheorem{defin}{\quad Definition}

\title{Minimum survival probabilities in a two-dimensional 
risk model perturbed by diffusion
}

\author{Chol-Ho Kim~~ and~~ Gwang-Ryong Han\\ \\
\small\textit{Department of Mathematics,} \\
\small\textit{Kim Il Sung University,  Pyongyang D.P.R.Korea}}

\date{July 20, 2012}

\begin{document}
\maketitle
%\hline

\begin{abstract}
In this paper we consider the finite time minimum survival probability and ultimate minimum survival probability in a two ? dimensional risk modal perturbed by diffusion
 Using some properties of the minimum survival probability we obtain the equation of the finite time minimum survival probability and ultimate minimum survival probability that they are satisfied and, the explicit expressions for ultimate minimum survival probability are given in a special case.

\bigskip

two-dimensional risk process
\end{abstract}

%\hline
%
%-----------------------1. Introduction-----------------
%
\section{Introduction}

In recent years a lot of work has been done with study of dependent classes of insurance business. Among other \cite{chan}(2003) considered most first, the ruin probability in the two-dimensional risk modal no perturbed by diffusion, and the 
problem to solve of sum ruin probability. And \cite{flo}(2008) studied the express of equation  satisfying  the ruin probability and  a method  to solve of the ruin probability in a special two- dimensional risk modal no perturbed by diffusion.
\cite{jun}(2007) considered the problem to obtain of Lundbergs lower bound for maximum ruin probability and studied the asymptotic estimation for finite time maximum ruin probability in a two-dimensional risk modal perturbed by diffusion. 
Consider the a two- dimensional risk process that is perturbed by diffusion:

\begin{eqnarray}
&& dR(t)=Cdt+\sigma dW(t)+\int_{|z|<1}\,\alpha (z)(\mu-\nu)(dt,dz) \nonumber\\
&& ~~~~~~~~~~~~ +\int_{|z|>1}\,\alpha (z) \mu (dt,dz) \\
&& R(0)=u \nonumber
\end{eqnarray}

\begin{equation*}
 C=(C_1,C_2)~~~~
 \mathbf{\sigma}=\left(
\begin{array}{ccc}

 \sigma _{11} & \sigma _{12} & \ldots  \,\sigma _{1d_1}\\
 \sigma _{21} & \sigma _{22} & \ldots  \,\sigma _{2d_1}
 \end{array}\right) 
\end{equation*} 
\begin{equation*}
 \mathbf{\alpha (Z)}=\left(
 \begin{array}{ccc}
 \alpha _{11}(Z) & \alpha _{12}(Z) & \ldots  \,\alpha _{1d_2}(Z)\\
 \alpha _{21}(Z) & \alpha _{22}(Z) & \ldots  \,\alpha _{2d_2}(Z)
 \end{array}\right)
 \end{equation*}
 where $u=(u_1,u_2)^{'}$ is the initial reserve of two insurance business and,$W(t)=(W_1(t),~\ldots \,,~W_{1d_1}(t))^{'}$ is independent  and $d_1-$ dimensional Wiener process vector, $\mu(dt,dz)=(\mu_1(dt,dz),\ldots\,,\mu_{d_2}(dt,dz))^{'}$ is independent and Poisson measure vector and that $E(\mu_i(dt,dZ))=\nu_i(dt,dZ)=\nu_i(dZ)dt$,and 
 $R(t)=(R_1(t),R_2(t))^{'}$ is the surplus reserve of two insurance business. Then process $R(t)=(R_1(t),R_2(t))^{'}$ 
 that is the surplus of two insurance company is homogeneous Markov process. \\
 We can rewriter $R(t)$
 \begin{eqnarray}
 && R(t)=u+Ct+dW(t)+\int_0^1 \int_{|Z|<1}\,\alpha(Z)(\mu - \nu)(ds,dZ)\nonumber\\
&& ~~~~~~~~~~~+\int_0^1 \int_{|Z|>1}\,\alpha(Z) \mu (ds,dZ)
 \end{eqnarray}
 Now one assume that $\int_{|z|>1}\,|\alpha_{k_i}(Z)| \nu_k(Z)<+ \infty $  for any $k\in \bar{1,d_2},\\i=1,2$\\
 Then $(1)$  can rewriter
 \begin{eqnarray*}
 && dR(t)=\ddot{C}(t)+\sigma dW(t)+\iint_{R_0^2}\,\alpha(Z)(\mu -\nu)(dt,dZ) ~~~~~~~~~~~~~~~~~~~~~~~~~~~~~~~~~~~~~~~~~~~~~~~~~~~~~~(1)^{'}
 \end{eqnarray*}
 where $\ddot{C}=C+\iint_{|z|\geq 1}\,\alpha(Z)\mu(dZ)$\\
 For the predicable three characteristics $(C,\sigma,\nu)$ of two-dimensional random process $\{R(t)\}_{t\in R,}$ and$~~\lambda=(\lambda_1,\lambda_2)^{'}\in R^2,$
 \begin{eqnarray*}
&& \textnormal{ln}\varphi(\lambda)=\{i(C,\lambda)+\frac{1}{2}(\lambda^{'}\sigma,\sigma^{'}\lambda)+\int\!\!\!\int_{R^2}\,[e^{i \textnormal{ln}(Z)}-1-\\
&& ~~~~~~~~~~~~~~~~~-i \lambda^{'} \alpha(Z)+\textnormal{I}\{|Z|<1\}]\nu(dZ)\}t 
 \end{eqnarray*}
We consider a  $\sigma-$ algebra $\Im_t$  that following way:\\
$\Im_t=\sigma\{R(s):s\leq t\}\bigcup N$\\
where $N=\sigma\{A_\alpha;P(A_\alpha)=0\}$
%-----------------------Definition1----------------- 
\begin{defin}
For $0\leq s,x\in R_+^2$
\begin{equation}
T_{s,x}=\textnormal{inf}\{t\geq s,R_1(t)<0 \textnormal{or} R_2(t)<0|R(s)=x\}
\end{equation}
is called minimum ruin time of state x at time s.
\end{defin}  
Then $T_{s,x}$ is $(\Im_t)-$ stop time.
%--------------------Definition2-----------------
\begin{defin}
For $0<s<T<+ \infty,x \in R_+^2$
\begin{eqnarray}
&& \phi_{\textnormal{min}}(s,x,T)=p(T_{s,x}<T)\nonumber\\
&& \Phi_{\textnormal{min}}{s,x,T}=1-\phi_{\textnormal{min}}(s,x,T)
\end{eqnarray}
is called respectively finite time minimum ruin probability and finite time minimum survival probability of state   at time s.
\end{defin}
And when $T=+\infty,$ i.e
\begin{eqnarray}
&& \phi_{\textnormal{min}}(s,x)=\lim_{T\rightarrow +\infty}\phi_{\textnormal{min}}(s,x;T)\nonumber\\
&& \Phi_{\textnormal{min}}(s,x)=\lim_{T\rightarrow +\infty}\Phi_{\textnormal{min}}(s,x;T)
\end{eqnarray}
is called respectively ultimate minimum ruin probability and ultimate minimum survival probability of state x at time s.
We leads to the equations that satisfied for finite time and ultimate minimum survival probability, and consider a method of solution for finite time minimum survival probability.
%----------2. Ultimate minimum survival probability equation.---------
\section{Ultimate minimum survival probability equation.}
%-----------------Theorem 1--------------------------
%
\begin{thr}
For any $0\leq s<t<T \phi_{\textnormal{min}}(t,R(t);T)$ is $(\Im_t-)$ martingale, and
\begin{equation*}
\Phi_{\textnormal{min}}(t,R(t);+\infty)=\Phi_{\textnormal{min}}(R(t))
\end{equation*}
\end{thr}
\textit{proof}
Let $A_{s,t}=\{\omega |\min_{s\leq t<T}\{R_1(t),R_2(t)\}>0\}$ ~~~~and
\begin{equation*}
I_{A_{s,T}}(\omega )=\left\{
\begin{array} {cc}
1 & \omega \in A_{s,t} \\
0 & \omega \notin A_{s,t}
\end{array}\right.
\end{equation*}
Then can write $\Phi_{\textnormal{min}}(s,x;T)=E\{I_{A_{s,T}}(\omega)|R(s)=x\}$\\
This can represent $\Phi_(s,R(s);T)=E\{I_{A_{s,T}}(\omega)|R(s)\}.$\\
For any $t(s<t<T)$ since $\{R(t)\}_t\in R^+$
\begin{eqnarray*}
&& E\{\Phi_{\textnormal{min}}(t,R(t);T)|\Im _t\}=E\{E\{I_{A_T}|R(t)\}|R(s)\}=E\{I_{A_T}|R(s)\}=\\
&& ~~~~~=\Phi_{\textnormal{min}}(s,R(s);T)
\end{eqnarray*}
i.e.$\Phi_{\textnormal{min}}(t,R(t);T)_{t \in [0,T)}$ is martingale , and taking into account that risk process  (2) is time homogeneous Markov process for $s<t$ ,we obtain that \\$p\{T_{s,x}<+ \infty \}=p\{T_{t,x}<+\infty \}$.\\
There fore $\Phi_{\textnormal{min}}(t,x;T)=p\{T_{s,x}\geq T\}=p\{T_{t,x}\geq T\}$ \\
Thus for $s,t (s<t) \\  \Phi_{\textnormal{min}}(t,x)=\Phi_{\textnormal{min}}(s,x)=\Phi_{\textnormal{min}}(x)$.
This implies theorem 1.\\
Now let 
$\Phi_{\textnormal{min}}(s,x;T)=\Phi_{\textnormal{min}}(s,x;),\;\; \Phi_{\textnormal{min}}(s,x;+\infty)=\Phi_{\textnormal{min}}(x)$.
%---------------------theorem2--------
\begin{thr}
If ~~~~$\Phi_{\textnormal{min}}(s,x)\in C^{1,2}\;\;([0,T]\times R_+^2)\;\;(\Phi_{\textnormal{min}}(x)\in C^{1,2}(R_+^2))$\\
for any $(s,x)\in [0,T]\times R_+^2\;\;(x\in R_+^2)$ then
\end{thr}
\begin{align}
& \left \{
 \begin{array}{l}
 \frac{\eth}{\eth s}\Phi_{\textnormal{min}}(s,x)+\tilde L \Phi_{\textnormal{min}}(s,x)=0 \\
 \Phi_{\textnormal{min}}(s,(x_1,0)^{'})=0,\Phi_{\textnormal{min}}(s,(x_2,0)^{'})=0,\Phi_{\textnormal{min}}(T,(x_1,x_2)^{'})=1
 \end{array}\right.
\\
& \left \{
 \begin{array}{l}
 \tilde L \Phi_{\textnormal{min}}(x)=0 \\
 \Phi_{\textnormal{min}}(x_1,0)=0,\Phi_{\textnormal{min}}(0,x_2)=0
 \end{array}\right.
 \end{align}
 where oparation $\tilde L$ is 
\begin{eqnarray}
&& \tilde L \Phi_{\textnormal{min}}(s,x)=\sum_{i=1}^2 C_i \frac{\eth \Phi_{\textnormal{min}}(s,x)}{\eth x_i}+\frac {1}{2}\sum_{i=1}^2 \sum_{j=1}^2 (\sum_{k=1}^{d_1}\sigma_{ik} \sigma_{jk})\cdot \nonumber\\
&& ~~~~~~~~\cdot  \frac{\eth ^2 \Phi_{\textnormal{min}}(s,x)}{\eth x_i \eth x_j}+\sum_{k=1}^{d_2} \iint _{R^2}\,(\Phi_{\textnormal{min}}(s,x)-\alpha_k(z))-\Phi_{\textnormal{min}}(s,x)-\nonumber\\
&& ~~~~~~~~-I_{|z|\leq1}\sum_{i=1}^2 \alpha_{ik}(z) \frac{\eth \Phi_{\textnormal{min}}(s,x)}{\eth x_i})\nu_k(dz)
\end{eqnarray}
\textit{proof}
~~~~~From theorem 1 since $\{\Phi_{\textnormal{min}}(t,R(t))\}_{t\in [0,T]}$ 
is martingale,applying ZhoZhuGyong transformation formula [1].
\begin{eqnarray}
&& \Phi_{\textnormal{min}}(t,R(t))=\Phi_{\textnormal{min}}(s,R(s))+\int_s^t \,\{\frac{\eth \Phi_{\textnormal{min}}(u,R(u))}{\eth u}+\sum_{i=1}^2 C_i \frac{\eth \Phi_{\textnormal{min}}(u,R(u))}{\eth x_i}\nonumber\\
&&~~~~+\frac {1}{2}\sum_{i=1}^2 \sum_{j=1}^2 \frac{\eth ^2 \Phi_{\textnormal{min}}(u,R(u))}{\eth x_i \eth x_j}\cdot \sum_{k=1}^{d_1}\sigma_{ik} \sigma_{jk}\}du+\int_s^t\,\sum_{i=1}^2 \frac{\eth\Phi_{\textnormal{min}}(u,R(u))}{\eth x_i}\cdot\nonumber\\
&&~~~~\cdot  \sum_{j=1}^{d_1} \sigma_{ij} d\omega_j(u)+\sum_{k=1}^{d_2} \int_s^t \,\iint_{|z|<1}\,[\Phi_{\textnormal{min}}(u,R(u)+\alpha_k(z) )-\nonumber\\
&& ~~~~-\Phi_{\textnormal{min}}(u,R(u))]-\sum_{i=1}^2 \alpha_{ik}(z)\frac{\eth}{\eth x_i}\Phi{\textnormal{min}}(u,R(u))\nu_k(du,dz)+\nonumber\\
&& ~~~~+\sum _{k=1}^{d_2}\int_s^t \,\iint_{|z|<1}\,[\Phi_{\textnormal{min}}(u,R(u)+\alpha_k(z) )-\Phi{\textnormal{min}}(u,R(u))](\mu_k- \nu_k )\nonumber\\
&& ~~~~(du,dz)+\sum _{k=1}^{d_2}\int_s^t \,\iint_{|z|<1}\,[\Phi_{\textnormal{min}}(u,R(u)+\alpha_k(z) )-\Phi{\textnormal{min}}(u,R(u))]\nonumber\\
&& (\mu_k- \nu_k )\mu_k (du,dz)
\end{eqnarray}
where 
\begin{equation*}
\mathbf{\alpha(z)}=\left(
\begin{array}{ccc}
\alpha_1(z)\\
\alpha_2(z)
\end{array}\right)
\end{equation*}
Thus from this theorem condition and result of theorem1, applying conditional Expectation of $\Im_s$, to (9), we obtain
\begin{eqnarray*}
&& E\{\Phi_{\textnormal{min}}(t,R(t))|\Im_s\}=\Phi_{\textnormal{min}}(s,R(s))+E\{ \int_s^t\,E\{ \frac{\eth}{\eth u}\Phi_{\textnormal{min}}(u,R(u))+\\
&& ~~~~+\sum_{i=1}^2 \frac {\eth}{\eth x_i} \Phi_{\textnormal{min}}(u,R(u))\cdot C_i+\frac {1}{2}\sum_{i=1}^2 \sum_{j=1}^2 (\sum_{k=1}^{d_1}\sigma_{ik} \sigma_{jk})\Phi_{\textnormal{min}}(u,R(u))+\\
&& ~~~~+ \iint_{R^2}\,\sum_{k=1}^{d_2} [ \Phi_{\textnormal{min}}(u,R(u)+\alpha_k(z) )-\Phi_{\textnormal{min}}(u,R(u))-I_{\{|z|\leq 1\}}\cdot \\
&&~~~~\cdot \sum_{i=1}^2\alpha_{ik}(z)\frac{\eth}{\eth x_i}\Phi_{\textnormal{min}}(u,R(u))]\nu_k(dz)\} du \}
\end{eqnarray*}
Therefore $\frac{\eth}{\eth u} \Phi_{\textnormal{min}}(u,x)+\tilde L \Phi_{\textnormal{min}}(u,x)=0$   , where $\tilde L$ is same to (8) and boundary conditions obtain from theorem 2 for finite time minimum survival probability.\\
In similar way applying same variation we obtain equation (7) for $\Phi_{\textnormal{min}}(R(t))$.\\                               
In case consider for Equation $(1)^{'}$. One obtain from result of theorem2. 
\begin{eqnarray}
&& \frac {\eth}{\eth s} \Phi_{\textnormal{min}}(s,x)+\sum_{i=1}^2 \frac {\eth}{\eth x_i} \Phi_{\textnormal{min}}(s,x)\cdot \tilde {C_i}+\frac {1}{2}\sum_{i=1}^2 \sum_{j=1}^2 (\sum_{k=1}^{d_1}\sigma_{ik} \sigma_{jk}) \cdot\nonumber\\
&& ~~~~\cdot \frac{\eth^2}{\eth x_i \eth x_j} \Phi_{\textnormal{min}}(s,x)+\sum_{k=1}^{d_2}\iint_{R^2}\,[\Phi_{\textnormal{min}}(s,x+\alpha_k(z) )-\nonumber\\
&& ~~~~-\Phi_{\textnormal{min}}(s,x)-\sum_{i=1}^2 \alpha_{ik}(z)\frac{\eth}{\eth x_i}\Phi_{\textnormal{min}}(s,x)]\nu_k(dz)=0
\end{eqnarray}
And taking into account KimZhuGyong [2]  that the random measure $\nu(dz)$ represent $\nu_k(dz)=\lambda_k p_k(z)$   for hyper density function $p_k(z)$   and $\alpha(z)=z$   in case $\alpha(z)=z$  we obtain from (10)
\begin{eqnarray*}
&& \frac {\eth}{\eth s} \Phi_{\textnormal{min}}(s,x)+\sum_{i=1}^2 \tilde C_i \frac {\eth}{\eth x_i} \Phi_{\textnormal{min}}(s,x)+\frac {1}{2}\sum_{i=1}^2 \sum_{j=1}^2 a_{ij}\frac{\eth^2 \Phi_{\textnormal{min}}(s,x)}{\eth x_i \eth x_j}-\\
&& ~~~~-a_0 \Phi_{\textnormal{min}}(s,x)+\iint_{R_+^2}\, \Phi_{\textnormal{min}}(s,z)\sum_{k=1}^d \lambda_k p_k(z-x)dz=0 
\end{eqnarray*}
where   
\begin{eqnarray*}
&&  \tilde C_i=C_i-\sum_{k=1}^{d_2} \lambda _k a_k, i=1,2\\
&& a_0= \sum _{k=1}^d \lambda _k, a_{ij}= \sum _{k=1}^{d_1} \sigma _{ik} \sigma _{jk}, i=1,2,j=1,2\\
&& a_k=\iint_{R^+}\,zp_k(z)dz, k=\bar {1,d_2}
\end{eqnarray*}

%-----------------------3. Finite time minimum survival probability-----------------
%
\section{Finite time minimum survival probability}

We consider two compound Poisson process
\begin{equation*}
x_1(t)=\sum_{k=1}^{M_1(t)} z_{1 k}, \; x_2(t)=\sum_{k=1}^{M_2(t)} z_{2 k}
\end{equation*}
where $\{z_{1k}\}, \{z_{2 k}\}$ are independent and have same distributions $F_1(z), F_2(z)$, respectively.
Then for two process we construct two-dimensional random measures in the following way:
\par
\textit{1. Case} $M_1(t)=M_2(t)=N(t)$ and ${N(t)}_{t \in R_k}$ is time homogeneous Poisson process that
$EN(t)=\lambda t$ for any $t \in R_+$ and $\Gamma \in B^2$ we define the measure $\mu$ that
\begin{equation}
\mu ((0, t], \Gamma )=\sum_{k=1}^{N(t)} I_\Gamma (z_{1k}z_{2k})
\end{equation}
Then measure $\mu$ have the following property

%
%-----------------Lemma 1--------------------------
%
\begin{lem}

\begin{eqnarray*}
&& (1) \; \mu ((0,t],\Gamma) \textnormal{is integral random measure on } B^2 \\
&& (2) \; \mu \textnormal{is non decreasing function and independent increment process of} t \\
&& (3) \; \nu ((0,t],\Gamma) = E(\mu((0,t],\Gamma)) = \lambda t \iint _ \Gamma d F_1(z_1) d F_2(z_2)
\end{eqnarray*}

\end{lem}

\par
\textit{2. Case} $M_1(t)=N_1(t)+N_3(t), M_2(t)=N_2(t)+N_3(t)$ and any $\{N_i(t)\}, \: i=1,2,3$ are independent
and have Poisson distribution with parameter $\lambda _i (i=1,2,3)$ we define the measures that

\begin{eqnarray}
&& \mu _1 ((0,t],\Gamma ^ {'}) = \sum_{k=1}^{N_1(t)} I_\Gamma (z_{1 k}, z_{3 k}) \nonumber\\
&& \mu _2 ((0,t],\Gamma ^ {'}) = \sum_{k=1}^{N_1(t)} I_\Gamma (z_{4 k}, z_{2 k}) \\
&& \mu _3 ((0,t],\Gamma ^ {'}) = \sum_{k=1}^{N_1(t)} I_\Gamma (z_{1 k}, z_{2 k}) \nonumber
\end{eqnarray}

Then the measures $\mu _ 1, \mu _2, \mu _ 3$ are independent and satisfy the results of Lemma 1.
Let $d_1 = 2, d_2 = 3$ in model (1) and (2), and using integral random measure (13) we can represent model
(1) in the following way:

\begin{equation}
\left \{
\begin{array} {lr}
dR(t) = Cdt+\sigma d W(t) - \iint_{R_+^2} z^{(1)} \mu_1 (dt,dz) -  \\
\qquad \qquad - \iint_{R_+^2} z^{(2)} \mu_2 (dt,dz) - \iint_{R_+^2} z \mu_3 (dt,dz) \\
R(0) = u
\end{array} \right.
\end{equation}

\noindent
where $R(t)=(R_1(t), R_2(t))^{'}, C=(C_1, C_2)^{'}$

\begin{eqnarray*}
&& \omega = \left(
\begin{array} {ccc}
\omega_1 (t) \\
\omega_2 (t)
\end{array}
\right)^{'} ~~~~~~~~~~~
\sigma = \left(
\begin{array} {ccc}
\sigma_1 & 0 \\
0 & \sigma_2
\end{array}
\right) \\
&& z^{(1)}=(z_1, 0) ^ {'}, z^{(2)} = (0, z_2) ^ {'}, z = (z_1, z_2)^{'} \\
&& u = (u_1, u_2) ^ {'}
\end{eqnarray*}

Then taking into account the equation of Theorem 2 which satisfied for model (14)
we obtain the following equation

\begin{eqnarray}
&& L \Phi_{min}(t,x) + \lambda_1 \int_{0}^{\infty}\Phi_{min}(t, x-z^{(1)})p_1(z_1)d z_1 + \nonumber\\
&& ~~~~~~~~~~~~~~ + \lambda_2 \int_{0}^{\infty}\Phi_{min}(t, x-z^{(2)})p_2(z_2)d z_2 + \nonumber\\
&& ~~~~~~~~~~~~~~ + \lambda_3 \int_{0}^{\infty}\!\!\!\int_{0}^{\infty}\Phi_{min}(t, x-z)p_1(z_1)p_2(z_2)d z_1 d z_2 = 0 \\
&& \Phi_{min} (t, (x_1, 0)^{'}) = \Phi_{min} (t, (0, x_2)^{'}) = 0, \Phi_{min} (T, (x_1, x_2)^{'}) = 1 \nonumber
\end{eqnarray}

\noindent
where operation $L$ is
\begin{eqnarray}
&& L= \frac{\partial}{\partial t} + \sum_{i=1}^{2}C_i \frac{\partial}{\partial x_i} + \frac{1}{2} \sum_{i=1}^{2} \partial_{i}^{2}
\frac{\partial ^ 2}{\partial_{x_i}^2}-\lambda \\
&& (\lambda = \lambda _1 + \lambda _2 + \lambda _3) \nonumber
\end{eqnarray}

Now we consider new part derivative operation
\begin{eqnarray}
&& L^*= -\frac{\partial}{\partial t} - \sum_{i=1}^{2}C_i \frac{\partial}{\partial x_i} + \frac{1}{2} \sum_{i=1}^{2} \partial_{i}^{2}
\frac{\partial ^ 2}{\partial_{x_i}^2}-\lambda \\
&& (\lambda = \lambda _1 + \lambda _2 + \lambda _3) \nonumber
\end{eqnarray}

%
%-----------------Theorem 3--------------------------
%
\begin{thr}

For $t(\tau \leq t \leq T_+^2), x \in R_+^2$, the solution of a part derivative equation
\begin{equation}
\left \{
\begin{array} {lr}
L^{*}k(t,x,z,\xi)=0 \\
k(t,(0,x_2)^{'})=0, ~~~~ k(t,(x,0)^{'})=0 \\
k(\tau, x;\tau, \xi)=\delta (\xi - x)
\end{array}
\right.
\end{equation}

\noindent
is obtained by the following formula :

\begin{eqnarray}
&& k(t,x,\tau, \xi)=\textnormal{exp}\left\{\beta (t-z) + <\alpha, x-\xi > \right\} \cdot \prod_{i=1}^{2} \frac{1}{\sqrt{2 \pi \sigma _i^2 (t-\tau)}} \nonumber\\
&& ~~~~~~~~~~~~~~~~~~ \left\{\textnormal{exp}\left\{-\frac{(x_i-\xi_i)^2}{2 \sigma_i^2 (t-\tau)}\right\}-
\textnormal{exp}\left\{-\frac{(x_i+\xi_i)^2}{2 \sigma_i^2 (t-\tau)}\right\}\right\}
\end{eqnarray}

\noindent
where $\alpha = (\alpha_1, \alpha_2)^{'}, \alpha_i = \frac{C_i}{\sigma_i^2},
\beta = -(\sum_{i=1}^{2}\frac{C_i}{\sigma_i^2}+\lambda), \xi = (\xi_1, \xi_2)^{'}$, and $<\cdot,\cdot>$ is the symbol
of scalar product.
\end{thr}

Denote the solution of (18) by $k(t,x;\tau, \xi)$. Taking into following transformation:
\begin{equation}
s(t,x;\tau, \xi)= e^{-<\alpha, x>}k(t,x;\tau, \xi)
\end{equation}
then the equation (18) can represent

\begin{align}
\frac{\partial s(t,x;\tau,\xi)}{\partial t}&-\left\{<C,\alpha>+\frac{1}{2}\alpha^{'}\sigma^{'}\alpha-\lambda\right\}s(t,x;\tau,\xi)+\nonumber\\
& \quad +\left[-C_1+\sigma_1^2\alpha_1\right]\frac{\partial s(t,x;\tau,\xi)}{\partial x_1}+
\left[-C_2+\sigma_2^2\alpha_2\right]\frac{\partial s(t,x;\tau,\xi)}{\partial x_2}+\\
& \quad +\frac{1}{2}\left[\sigma_1^2\frac{\partial^2s(t,x;\tau,\xi)}{\partial x_1^2}+
\sigma_2^2\frac{\partial^2s(t,x;\tau,\xi)}{\partial x_2^2}\right]=0\nonumber
\end{align}

Denoting $\alpha_i=\frac{C_i}{\sigma_i^2}, \quad i=1,2,
\quad \beta=-\left(\lambda+\frac{C_1^2}{2\sigma_1^2}+\frac{C_2^2}{2\sigma_2^2}\right)$, then for $t(\tau \leq t \leq T)$ and
$x,\xi \in R_+^2$ the equation (21) can represent the following way:

\begin{equation}
\left \{
\begin{array} {lr}
\left\{-\frac{\partial}{\partial t}-\frac{1}{2}\left[\sigma_1^2\frac{\partial^2}{\partial x_1^2}+
\sigma_2^2\frac{\partial^2}{\partial x_2^2}\right]+\beta\right\}s(t,0;\tau,\xi)=0 \\
s(t,x;\tau,\xi)=0 \\
s(\tau,x;\tau,\xi)=e^{-<\alpha,x>}\delta(x-\xi)
\end{array} \right.
\end{equation}

If $\Gamma (t,x;\tau,\xi)$ is an extended solution of (21) in definite domain $[\tau, T]\times R^2$, we can represent

\begin{equation}
\Gamma (t,x;\tau,\xi) = 
\left \{
\begin{array} {lr}
s(t,(x_1,x_2);\tau,\xi), & x_1 \geq 0, x_2 \geq 0 \\
-s(t,(-x_1,x_2);\tau,\xi), & x_1>0, x_2 \geq 0 \\
-s(t,(x_1,-x_2);\tau,\xi), & x_1 \geq 0, x_2<0 \\
s(t,(-x_1,x_2);\tau,\xi), & x_1<0, x_2<0
\end{array} \right.
\end{equation}

Then initial condition of (22) can represent.

\begin{equation*}
\Gamma(\tau,x;\tau,\xi) =
\left \{
\begin{array} {lr}
e^{<\alpha,x>}\delta(x_1-\xi_1,x_2-\xi_2), & x_1 \geq 0, x_2 \geq 0 \\
-e^{\alpha_1 x_1 - \alpha_2 x_2}\delta(x_1-\xi_1,x_2+\xi_2), & x_1>0, x_2<0 \\
-e^{-\alpha_1 x_1 + \alpha_2 x_2}\delta(x_1+\xi_1,x_2-\xi_2), & x_1<0, x_2>0 \\
e^{<\alpha,x>}\delta(x_1+\xi_1,x_2+\xi_2), & x_1<0, x_2<0
\end{array} \right.
\end{equation*}

Thus the boundary value problem of (21) leads to the initial value problem of $\Gamma (t,x;\tau,\xi)$

\begin{equation}
\left \{
\begin{array} {lr}
\left[-\frac{\partial}{\partial t}+\frac{1}{2}\left(\sigma_1^2\frac{\partial^2}{\partial x_1^2}+\sigma_2^2\frac{\partial^2}{\partial x_2^2}\right)
+\beta\right]\Gamma(t,x;z,\xi)=0 \\
\Gamma(t,x;z,\xi)=0=\varphi(x,\xi)
\end{array} \right.
\end{equation}

From the solving method of Kolmogorov-Feller equation

\begin{align*}
\Gamma (t,x;z,\xi) &= \varphi (x,\xi)\ast \left( \frac{1}{2 \pi} \right) ^2 \iint_{R^2} e^{-i<\theta, X>}\cdot \\
& \qquad \quad \cdot \textnormal{exp} \left\{-\frac{1}{2}\left(\sigma_1^2 \theta_1^2 + \sigma_2^2 \theta_2^2 \right) (t-\tau) +
\beta(t-\tau)\right\}d\theta_1 \theta_2 \\
&=  \varphi (x, \xi) \ast \prod_{i=1}^{2}\frac{1}{\sqrt{2 \pi \sigma_i^2(t-\tau)}}\textnormal{exp}
\left\{-\frac{x_i^2}{2\sigma_i^2(t-\tau)}\right\}\textnormal{exp}\{\beta(t-\tau)\}
\end{align*}
where $\ast$ is convolution symbol and $\theta = (\theta_1, \theta_2)^{'}$.
Thus,
\begin{align*}
\Gamma (t,x;z,\xi) &= \iint_{R^2}\varphi (x,y,\xi) \prod_{i=1}^{2}\frac{1}{\sqrt{2 \pi \sigma_i(t-\tau)}} \cdot \\
& \qquad \quad \cdot \textnormal{exp}\left\{-\frac{y_I^2}{2\sigma_i(t-\tau)}\right\}\textnormal{exp}\{\beta(t-\tau)\}d y_1 y_2= \\
& =\textnormal{exp}\{-<\alpha,\xi>\}\frac{1}{\sqrt{2 \pi (t-\tau) \sigma_1 \sigma_2}} \cdot \\
& \qquad \quad \cdot \{\textnormal{exp}\left\{-\frac{\frac{(x_1-\xi_1)^2}{\sigma_1^2}+\frac{(x_2-\xi_2)^2}{\sigma_2^2}}{2(t-\tau)} \right\}-\\
& \qquad \quad - \textnormal{exp}\left\{-\frac{\frac{(x_1-\xi_1)^2}{\sigma_1^2}+\frac{(x_2+\xi_2)^2}{\sigma_2^2}}{2(t-\tau)} \right\}- \\
& \qquad \quad - \left\{-\frac{\frac{(x_1-\xi_1)^2}{\sigma_1^2}+\frac{(x_2-\xi_1)^2}{\sigma_2^2}}{2(t-\tau)} \right\} + \\
& \qquad \quad + \left\{-\frac{\frac{(x_1-\xi_1)^2}{\sigma_1^2}+\frac{(x_2+\xi_2)^2}{\sigma_2^2}}{2(t-\tau)} \right\} \}
\textnormal{exp}\{\beta (t-\tau)\}
\end{align*}

Taking into account he (20) and (23), the solution of equation (18) is
\begin{align*}
K(t,x;z,\xi) &= e^{<\alpha,X>}S(t,x;z,\xi) =\\
&= \textnormal{exp}\{<x,x-\xi>+\beta (t-\tau)\}\frac{1}{2 \pi (t-\tau) \sigma_1 \sigma_2} \cdot \\
& \qquad \quad \cdot \prod_{i=1}^2\left\{\textnormal{exp}\left\{-\frac{(x_i-\xi_i)^2}{2(t-\tau)\sigma_i^2}\right\}-
\textnormal{exp}\left\{-\frac{(x_i+\xi_i)^2}{2(t-\tau)\sigma_i^2}\right\}\right\}
\end{align*}
%
%-----------------Theorem 4--------------------------
%
\begin{thr}

A part derivative equation (15) is equivalent to the following integral equation
\begin{equation}
\Phi_{min}(\tau,\xi)=\int_\tau^T\iint_{R_+^2}\Phi_{min}(t,x)G(t,x;\tau,\xi)dxdt+F(\tau,\xi)
\end{equation}
\noindent
where
\begin{align*}
F(\tau,\xi)&=\iint_{R_+^2}K(T,x;\tau,\xi)d x \\
G(t,x;\tau,\xi) &= -\lambda_1 \int_{x_1}^{+\infty}K(t,(z_1,z_2);\tau,\xi)P_1(x_1-z_1)d z_1- \\
& -\lambda_2 \int_{x_2}^{+\infty}K(t,(x_1,z_2);\tau,\xi)P_2(x_2-z_2)d z_2- \\
& -\lambda_3 \int_{x_1}^{+\infty} \int_{x_2}^{+\infty}K(t,(z_1,z_2);\tau,\xi)P_1(x_1-z_1)P_2(x_2-z_2)d z_1 d z_2
\end{align*}

\end{thr}

%
%-----------------Theorem 5--------------------------
%
\begin{thr}
Integral equation (25) in domain $\Sigma = \left\{(t,\xi) \in [\tau,t]\times R_+^2 \right\}$ has a unique solution
\begin{equation}
\Phi_{min}(\tau,\xi)=F(\tau,\xi)+\sum_{k=1}{\infty}\int_\tau^T\iint_{R_+^2}G^\infty(t,x;z,\xi)F(t,x)d x d t
\end{equation}
\noindent
where
\begin{align*}
G^{(1)}(t,x;z,\xi)&=G(t,x;z,\xi) \\
G^{(k)}(t,x;z,\xi)&=\int_\tau^T\iint_{R_+^2}G^{k-1}(t,x;s,z)d z d s
\end{align*}

\end{thr}

\end{document}